\documentclass[a4paper,12pt]{amsart}
\usepackage{amssymb}
\usepackage{ifthen}
\usepackage{graphicx}
\usepackage{multicol}        
\usepackage{epsfig}
\nonstopmode \numberwithin{equation}{section}
\newtheorem{thm}[equation]{Theorem}
\newtheorem{cor}[equation]{Corollary}
\newtheorem{lem}[equation]{Lemma}
\newtheorem{prop}[equation]{Proposition}

\newtheorem{conj}{Conjecture}

\theoremstyle{definition}

\newtheorem{prob}[equation]{Problem}

\newenvironment{rem}{%
\bigskip
\noindent \textsl{{\sl Remark. }}}{\bigskip}
\newenvironment{rems}{%
\bigskip
\noindent \textsl{{\sl Remarks. }}}{\bigskip}

\newcounter {own}
\def\theown {\thesection       .\arabic{own}}

\newenvironment{pf}[1][]{%
 \vskip 3mm
 \noindent
 \ifthenelse{\equal{#1}{}}%
  {{\slshape Proof. }}%
  {{\slshape #1.} }%
 }%
{\qed\bigskip}

\newcounter{alphabet}
\newcounter{tmp}
\newenvironment{Thm}[1][]{\refstepcounter{alphabet}%
\bigskip%
\noindent%
{\bf Theorem \Alph{alphabet}}%
\ifthenelse{\equal{#1}{}}{}{ (#1)}%
{\bf .} \itshape}{\vskip 8pt}

\newcommand{\ID}{{\mathbb D}}

\newcommand{\IC}{{\mathbb C}}

\def\be{\begin{equation}}
\def\ee{\end{equation}}

\newcommand{\bee}{\begin{enumerate}}
\newcommand{\eee}{\end{enumerate}}

\newcommand{\blem}{\begin{lem}}
\newcommand{\elem}{\end{lem}}
\newcommand{\bthm}{\begin{thm}}
\newcommand{\ethm}{\end{thm}}
\newcommand{\bcor}{\begin{cor}}
\newcommand{\ecor}{\end{cor}}
\newcommand{\beg}{\begin{examp}}
\newcommand{\eeg}{\end{examp}}
\newcommand{\begs}{\begin{examples}}
\newcommand{\eegs}{\end{examples}}
\newcommand{\bdefe}{\begin{defin}}
\newcommand{\edefe}{\end{defin}}
\newcommand{\bprob}{\begin{prob}}
\newcommand{\eprob}{\end{prob}}
\newcommand{\bei}{\begin{itemize}}
\newcommand{\eei}{\end{itemize}}

\newcommand{\bcon}{\begin{conj}}
\newcommand{\econ}{\end{conj}}
\newcommand{\bcons}{\begin{conjs}}
\newcommand{\econs}{\end{conjs}}
\newcommand{\bprop}{\begin{prop}}
\newcommand{\eprop}{\end{prop}}
\newcommand{\br}{\begin{rem}}
\newcommand{\er}{\end{rem}}
\newcommand{\brs}{\begin{rems}}
\newcommand{\ers}{\end{rems}}
\newcommand{\bo}{\begin{obser}}
\newcommand{\eo}{\end{obser}}
\newcommand{\bos}{\begin{obsers}}
\newcommand{\eos}{\end{obsers}}
\newcommand{\bpf}{\begin{pf}}
\newcommand{\epf}{\end{pf}}
\newcommand{\ba}{\begin{array}}
\newcommand{\ea}{\end{array}}
\newcommand{\beq}{\begin{eqnarray}}
\newcommand{\beqq}{\begin{eqnarray*}}
\newcommand{\eeq}{\end{eqnarray}}
\newcommand{\eeqq}{\end{eqnarray*}}

\newcommand{\lra}{\longrightarrow}

\newcommand{\ds}{\displaystyle}

\newcounter{minutes}
\divide\time by 60
\newcounter{hours}
\multiply\time by 60 \addtocounter{minutes}{-\time}

\begin{document}
\bibliographystyle{amsplain}
\title{Region of variability for
exponentially convex univalent functions}
\author{S. Ponnusamy}
\address{S. Ponnusamy, Department of Mathematics,
Indian Institute of Technology Madras, Chennai-600 036, India.}
\email{samy@iitm.ac.in}
\author{A. Vasudevarao}
\address{A. Vasudevarao, Department of Mathematics,
Indian Institute of Technology Madras, Chennai-600 036, India.}
\email{alluvasudevarao@gmail.com}
\author{M. Vuorinen}
\address{M. Vuorinen, Department of Mathematics,
FIN-20014 University of Turku, Finland.}
\email{vuorinen@utu.fi}

\date{
\texttt{File:~\jobname .tex,
          printed: \number\year-\number\month-\number\day,
          \thehours.\ifnum\theminutes<10{0}\fi\theminutes}
}
\begin{abstract}
For $\alpha\in\IC\setminus \{0\}$  let $\mathcal{E}(\alpha)$
denote the class of all univalent functions $f$ in the unit disk $\mathbb{D}$ and is given by
$f(z)=z+a_2z^2+a_3z^3+\cdots$, satisfying
$$ {\rm Re\,} \left (1+ \frac{zf''(z)}{f'(z)}+\alpha zf'(z)\right )>0
\quad \mbox{ in  ${\mathbb D}$}.
$$
For any fixed $z_0$ in the unit disk $\mathbb{D}$ and
$\lambda\in\overline{\mathbb{D}}$,
we determine the region of variability $V(z_0,\lambda)$ for
$\log f'(z_0)+\alpha f(z_0)$ when $f$ ranges over the class
$$\mathcal{F}_{\alpha}(\lambda)=\left\{f\in\mathcal{E}(\alpha) \colon
f''(0)=2\lambda-\alpha
\right\}.
$$
We geometrically illustrate the region of variability $V(z_0,\lambda)$
for several sets of parameters using Mathematica.
In the final section of this article we propose some open problems.
\end{abstract}

\thanks{The first two authors were supported by NBHM (DAE, sanction No. 48/2/2006/R\&D-II).
The work was completed during the visit of the first two authors to the University of Turku, Finland. }

\maketitle
\pagestyle{myheadings}
\markboth{S. Ponnusamy, A. Vasudevarao and M. Vuorinen}{Regions of variability}

\section{Introduction}\label{sec1}

Let $\ID=\{z\in\mathbb{C}\colon |z|<1\}$ be the unit disk in the complex plane
$\IC$. We denote the class of analytic functions in $\ID$ by
${\mathcal H}$ which we think of as a topological
vector space endowed with the topology of uniform convergence over
compact subsets of $\mathbb{D}$. Let $\mathcal{A}$ denote the family
of functions $f$ in ${\mathcal H}$ normalized by $f(0)=f'(0)-1=0$.
A function $f\in\mathcal{A}$ is said to be in the family $\mathcal{S}$
if it is univalent in $\ID$. Denote by $ {\mathcal S}^*$ the subclass of functions
$\phi\in{\mathcal A}$ such that $\phi$ maps
$\mathbb{D}$ univalently onto a domain $\Omega=\phi(\mathbb{D})$
that is starlike with respect to the origin. That is,
$t\phi(z)\in\phi(\mathbb{D})$ for each  $t\in [0,1]$.
It is well known that $\phi\in {\mathcal S}^*$
if and only if
$${\rm Re\,}\left (\frac{z\phi '(z)}{\phi (z)}\right )>0, \quad z\in \ID.
$$
Functions in ${\mathcal S}^*$ are referred to as starlike functions.
It is well known that $\phi\in\mathcal{A}$ maps $\mathbb{D}$
univalently  onto a convex domain, denoted by $\phi \in {\mathcal
C}$, if and only if $z\phi'\in{\mathcal S}^*$. Functions in
${\mathcal C}$ are referred to as normalized convex functions. We refer to the
books \cite{Du,Goodman's_book,Pommerenke-book2} for a detailed discussion on these
two classes.

A function $f\in \mathcal{H}\,$ is called exponentially convex if $f$ is univalent in $\ID$  and  $e^{f(z)}$
maps $\ID$ onto a convex domain.  For $\alpha\in\IC\setminus\{0\}$, the
family $\mathcal{E}(\alpha)$ of $\alpha$-exponential functions was
introduced in \cite{Arango-Mejia-Rusch-97}. A function $f\in\mathcal{S}$ is said
to be in $\mathcal{E}(\alpha)$ if $F(\ID)$ is a convex domain, where
$F(z)=e^{\alpha f(z)}$.
Although Arango {\it et al} \cite{Arango-Mejia-Rusch-97}  introduced and studied exponentially convex functions
in 1997,  no attempt has been made to study many other properties of the class $\mathcal{E}(\alpha)$ until the
present article. In this article, we initiate certain issues related to $\mathcal{E}(\alpha)$. We
now recall a number of basic properties of the class
$\mathcal{E}(\alpha)$ from \cite{Arango-Mejia-Rusch-97}.

\begin{Thm}{\rm (\cite[Theorem 1]{Arango-Mejia-Rusch-97})}
Let $\alpha\in\IC\setminus\{0\}$. A function $f$ is in $\mathcal{E}(\alpha)$ if and only if
${\rm Re\,} P_f(z)>0$ in $\ID$, where
\be\label{pvdev-9-eq1}
P_f(z)=1+\frac{zf''(z)}{f'(z)}+\alpha zf'(z).
\ee
\end{Thm}

\begin{Thm}{\rm (\cite[Theorem 2]{Arango-Mejia-Rusch-97})}
Let $f\in \mathcal{E}(\alpha)$. Then $f(\ID)$ is convex in the $\overline{\alpha}-$
direction (and therefore close-to-convex). It is not necessarily starlike univalent.
\end{Thm}

Because each $f\in\mathcal{E}(\alpha)$ can be represented in the
form
$$e^{\alpha f(z)}=1+\alpha g(z), \quad g\in\mathcal {A},
$$
as shown in \cite{Arango-Mejia-Rusch-97}, $g$ belongs to the class ${\mathcal C}(\alpha)$
of (normalized) convex univalent functions with $-1/\alpha\notin g(\ID)$. Here ${\mathcal C}(\alpha)$ denotes
the class of convex functions of order $\alpha$.

\begin{Thm}{\rm (\cite[Theorem 3]{Arango-Mejia-Rusch-97})}
For $\alpha\in\IC\setminus \{0\}$ we have
$$
\mathcal{E}(\alpha)=\left\{\frac{1}{\alpha}\log(1+\alpha g) \colon
g\in\mathcal{C}(\alpha)\right\}.
$$
\end{Thm}

In \cite[Theorem 4]{Arango-Mejia-Rusch-97},  it is also observed that
$\mathcal{E}(\alpha)=\emptyset$ if $|\alpha|>2$ and for $|\alpha|=2$,
$\mathcal{E}(\alpha)$ consists of the functions
$$
f(z)=\frac{1}{\alpha}\log \frac{2+\alpha z}{2-\alpha z}.
$$
Throughout the discussion, we assume that $|\alpha|\leq 2$.

For $f\in{\mathcal F}_\alpha $, we denote by $\log f'$ the
single-valued branch of the logarithm of $f'$ with $\log f'(0)=0$.
The Herglotz formula for analytic functions that have positive real
part in the unit disk shows that if $f\in\mathcal{E}(\alpha)$, then
there exists a unique positive unit measure $\mu$ on $(-\pi,\pi]$
such that
$$1+\frac{zf''(z)}{f'(z)}+ \alpha zf'(z)=
\int_{-\pi}^{\pi}\frac{1+ze^{-it}}{1-ze^{-it}}\,d\mu(t), \quad z\in \ID.
$$
A computation gives that
$$\log f'(z)+\alpha f(z)=
\int_{-\pi}^{\pi}\log\left(\frac{1}{1-ze^{-it}}\right)^2 d\mu(t),
$$
or equivalently we can write
$$
\log f'(z)+\alpha f(z)=2\int_0^1\frac{\omega(tz)}{1-\omega(tz)}\,\frac{dt}{t}
$$
for some $\omega\in{\mathcal B}_0$.
Here ${\mathcal B}_0$ denotes the class of analytic functions $\omega$
in $\mathbb{D}$ such that $| \omega(z)|\leq 1$ in $\mathbb{D}$ and
$\omega(0)=0$. Consequently, for each $f\in\mathcal{E}(\alpha)$ there
exists an $\omega_f \in {\mathcal B}_0$ of the form
\be\label{pvdev-9-eq2}
\omega_f(z)=\frac{P_f(z)-1}{P_f(z)+1},\quad z\in\mathbb{D},
\ee
and conversely. It is a simple exercise to  see that
\be\label{pvdev-9-eq2a}
P'_f(0)=2\omega'_f(0)=f''(0)+\alpha.
\ee
Suppose that $f\in\mathcal{E}(\alpha)$. Then, a simple application of the Schwarz lemma
(see for example \cite{Du,samy-book3,samy-herb}) shows that
$$|P'_f(0)|= |f''(0)+\alpha|\leq 2,
$$
because $|\omega'_f(0)|\leq 1$. That is $f''(0)=2\lambda-\alpha$ for some
$\lambda\in\mathbb{\overline{D}}=\{z\in\mathbb{C}\colon |z|\leq 1\}$.
For $\omega_f\in{\mathcal B}_0$, we define the function
$g:\mathbb{D}\lra\overline{\mathbb{D}}$ by

\be\label{pvdev-9-eq2b}
g(z) =
\left\{
\ba{ll}
\ds \frac{\frac{\omega_f(z)}{z}-\lambda}{1-\overline{\lambda}\frac{\omega_f(z)}{z}}
&\mbox{if $|\lambda|<1$}, \\[6mm]
\ds 0
& \mbox{if $|\lambda|=1$} .
\ea \right.
\ee
Then we see that
\be\label{pvdev-9-eq2d}
g'(0) =
\left\{
\ba{ll}
\ds \frac{\omega_f''(0)}{2(1-|\lambda|^2)}
&\mbox{if $|\lambda|<1$}, \\[6mm]
\ds 0
& \mbox{if $|\lambda|=1$} .
\ea \right.
\ee
From  (\ref{pvdev-9-eq2}) and (\ref{pvdev-9-eq2a}), we find that
\begin{eqnarray}\label{pvdev-9-eq2e}
2\omega''_f(0)+(P'_f(0))^2=P''_f(0).
\end{eqnarray}
As  $f''(0)=2\lambda-\alpha$, using (\ref{pvdev-9-eq1}) and (\ref{pvdev-9-eq2e})
we obtain
\be\label{pvdev-9-eq2c}
\omega''_f(0)=f'''(0)-6\lambda(\lambda-\alpha)-2\alpha^2.
\ee
From (\ref{pvdev-9-eq2d}) and (\ref{pvdev-9-eq2c}) we  also note that $|g'(0)|\leq 1$ if and only if
\begin{equation*}
f'''(0)=2[(1-|\lambda|^2)a+ 3\lambda(\lambda-\alpha)+\alpha^2]
\end{equation*}
for some $a\in\mathbb{\overline{D}}$.
Consequently, for
$\lambda\in\overline{\mathbb{D}}$ and for
$z_0\in\mathbb{D}$ fixed, it is natural to introduce
\beqq
{\mathcal F}_\alpha(\lambda) & = &
\left\{f\in\mathcal{E}(\alpha) \colon  f''(0)=2\lambda-\alpha\right\},
\eeqq
and
$$
V (z_0,\lambda) = \{\log f'(z_0)+\alpha f(z_0)\colon f\in {\mathcal F}_\alpha(\lambda)\}.
$$
From (\ref{pvdev-9-eq2a}) and the normalization condition
in the class ${\mathcal F}_\alpha(\lambda)$, we observe that
$\omega'_f(0)=\lambda$.

In the recent past,  several authors
(see \cite{samy-vasudev-vuorinen2,Yanagihara2} and references there in)
have studied region of variability problems for several subclasses of ${\mathcal S}$. The main aim of this article is  to determine the region of variability $V (z_0,\lambda)$ of $\log f'(z_0)+\alpha
f(z_0)$ when $f$ ranges over the class ${\mathcal F}_\alpha(\lambda)$.
Our main theorem here is Theorem \ref{pvdev-9-th1} and at the end we graphically illustrate the region of
variability for several sets of parameters. Finally,
we propose some open problems on exponentially convex univalent functions.

\section{The Basic properties of $V (z_0,\lambda)$ and  the  Main result}\label{sec2}

To state our main theorem, we need some  preparation. For a
positive integer $p$, let
$$({\mathcal S}^*)^p=\{f=f_0^p\colon f_0\in {\mathcal S}^* \}
$$
and recall the following result from \cite{Yanagihara2}.

\blem\label{pvdev-9-lem1} Let $f$ be an analytic function in
${\mathbb D}$ with $f(z) = z^p + \cdots $. If
$$ {\rm Re} \,  \left( 1+ z \frac{f''(z)}{f'(z)} \right)> 0 , \quad z \in {\mathbb D} ,
$$
then $f \in ({\mathcal S}^*)^p$.
\elem

Now, we list down some basic properties of $V (z_0,\lambda)$.
\bprop\label{pvdev-9-pro3}
We have
\bee
\item $V (z_0,\lambda)$ is a compact subsect of $\mathbb{C}$.
\item $V (z_0,\lambda)$ is a convex subsect of $\mathbb{C}$.
\item For $|\lambda|=1$ or $z_0=0$,
\be\label{pvdev-9-eq20}
V (z_0,\lambda)=\left \{-2\log(1-\lambda z_0)\right \}.
\ee
\item For $|\lambda|<1$ and $z_0\neq 0$, $V (z_0,\lambda)$ has
$-2\log(1-\lambda z_0)$ as an interior point.
\eee
\eprop
\bpf (1) Since ${\mathcal F}_\alpha (\lambda)$ is a
compact subset of $\mathcal{H}$, it follows that
$V (z_0,\lambda)$ is also a compact subset of $\mathbb{C}$.

(2) If $f_0,f_1\in {\mathcal F}_\alpha (\lambda)$ and $0\leq t \leq 1$, then
the function $f_t$ satisfying 
$$
\log f'_t(z)+\alpha f_t(z)=(1-t)(\log f'_0(z)
+\alpha f_0(z))+t(\log f'_1(z)+\alpha f_1(z))
$$
is evidently in ${\mathcal F}_\alpha(\lambda)$. Also, because of the
representation of $f_t$, we see easily that  the set $V (z_0,\lambda)$ is a
convex subset of $\mathbb{C}$.

(3) If $z_0=0$, then (\ref{pvdev-9-eq20}) trivially holds. If $|\lambda|=|\omega'_f(0)|=1$,
then it follows from the Schwarz lemma (see for example \cite{Du,samy-book3,samy-herb})
that $\omega_f(z)=\lambda z$, which implies
$$P_f(z)=\frac{1+\lambda z}{1-\lambda z}.
$$
A simple computation yields
$$
\log f'(z)+\alpha f(z)=-2\log(1-\lambda z).
$$
Consequently,
$$V (z_0,\lambda)
=\left \{-2\log(1-\lambda z_0)\right \}.
$$

(4)  For $\lambda\in\mathbb{D}$, $ z_0\in\mathbb{D}\setminus\{0\}$,
and $a\in\overline{\mathbb{D}}$, we define
$$\delta(z,\lambda) = \frac{z+\lambda}{1+\overline{\lambda}z}.
$$
Recall that $|g'(0)|\leq 1$ if and only if
$f'''(0)=2[(1-|\lambda|^2)a+ 3\lambda(\lambda-\alpha)+\alpha^2]$ for some $a\in\mathbb{\overline{D}}$.
Now by applying the Schwarz lemma to $g(z)$ defined in (\ref{pvdev-9-eq2b}) we obtain
\begin{equation}\label{pvdev-9-eq20a}
g(z)=az\quad\mbox{ for some } |a|=1.
\end{equation}
A  simple computation of (\ref{pvdev-9-eq20a}) gives
\be\label{pvdev-9-eq4}
\log H'_{a,\lambda}(z) +\alpha H_{a,\lambda}(z) =
\int_0^z\frac{2\delta(a\zeta, \lambda)}
{1-\delta(a\zeta, \lambda)\zeta}\,d\zeta ,\quad
z\in\mathbb{D}.
\ee

First we claim that  $H_{a,\lambda}$ satisfying (\ref{pvdev-9-eq4}) belongs ${\mathcal F}_\alpha(\lambda)$.
For this, we compute
\begin{eqnarray*}
1+\frac{zH''_{a,\lambda}(z)}{H'_{a,\lambda}(z)}+\alpha  H'_{a,\lambda}(z)
& = & \frac{1+\delta(az,\lambda)z}{1-\delta(az,\lambda)z}.
\end{eqnarray*}
As $\delta(az,\lambda)$ lies in the unit disk $\mathbb{D}$,
$H_{a,\lambda}\in {\mathcal F}_\alpha(\lambda)$ and the claim follows.
Also we observe that
\be\label{pvdev-9-eq5}
\omega_{H_{a,\lambda}}(z)=z\delta(az,\lambda).
\ee

Next we claim that the mapping
${\mathbb D} \ni a\mapsto \log H'_{a,\lambda}(z_0)+\alpha H_{a,\lambda}(z_0) $
is a non-constant analytic function of $a$
for each fixed $z_0 \in {\mathbb D} \backslash \{ 0 \}$ and
$\lambda\in\mathbb{D}$. To do this, we put
\begin{eqnarray*}
h(z) & = & \left .\frac{1}{(1-|\lambda|^2)}
\frac{\partial}{\partial a}\left\{\frac{}{}\log H'_{a,\lambda} (z)
+\alpha H_{a,\lambda}(z)\right\}\right |_{a=0}.
\end{eqnarray*}
A computation gives
\begin{eqnarray*}
h(z) & = &  2\int_0^z \frac{\zeta}{(1-\lambda\zeta)^2}\,d\zeta=z^2+\cdots
\end{eqnarray*}
from which it is easy to see that
$${\rm Re} \,\left\{1+\frac{zh''(z)}{h'(z)}\right\}=
{\rm Re} \,\left\{\frac{2}{1-\lambda z}\right\}> 1,
\quad z\in\mathbb{D}.
$$
By Lemma \ref{pvdev-9-lem1} there exists a function $h_0\in
{\mathcal S}^*$  with $h=h_0^2$. The univalence of $h_0$ together with
the condition $h_0(0)=0$ implies that $h(z_0)\neq 0$ for $z_0 \in {\mathbb
D}\setminus \{0\}$. Consequently, the mapping
${\mathbb D} \ni a\mapsto \log H'_{a,\lambda}(z_0) +\alpha H_{a,\lambda}(z_0)$
is a non-constant analytic function of $a$ and hence, it is an open mapping. Thus,
$V (z_0,\lambda)$ contains the open set
$\{\log H'_{a,\lambda}(z_0)+\alpha H_{a,\lambda}(z_0)\colon |a|<1\}$.
In particular,
$$\log H'_{0,\lambda}(z_0)+\alpha H_{0,\lambda}(z_0) = -2\log (1-\lambda z_0)
$$
is an interior point of
$\{\log H'_{a,\lambda}(z_0)+\alpha H_{a,\lambda}(z_0)\colon
a\in\mathbb{D}\}\subset V (z_0,\lambda)$.
\epf

We remark that, since $V (z_0,\lambda)$ is a compact convex subset of
$\mathbb{C}$ and has nonempty interior, the boundary
$\partial{V (z_0,\lambda)}$ is a Jordan curve and $V (z_0,\lambda)$
is the union of $\partial{V (z_0,\lambda)}$ and its inner domain.
Now we state our main result and the proof will be presented in
Section \ref{sec3}.

\bthm \label{pvdev-9-th1} For $\lambda\in\mathbb{D}$, $\alpha\in\IC\setminus \{0\}$ and
$z_0\in\mathbb{D}\setminus \{0\}$, the boundary
$\partial{V (z_0,\lambda)}$ is the Jordan curve given by
\begin{eqnarray*}
(-\pi,\pi]\ni \theta \mapsto \log H'_{e^{i\theta},\lambda}(z_0)+
\alpha H_{e^{i\theta},\lambda}(z_0) & =&
\int_0^{z_0}\frac{2\delta(e^{i\theta}\zeta, \lambda)}
{1-\delta(e^{i\theta}\zeta, \lambda)\zeta} \,d\zeta.
\end{eqnarray*}
If  $ f(z_0)= H_{e^{i\theta},\lambda}(z_0)$ for some
$f\in{\mathcal F}_\alpha(\lambda)$ and $\theta\in (-\pi,\pi]$, then
$f(z)=H_{e^{i\theta},\lambda}(z)$.
\ethm

\section{Preparation for the proof of Theorem \ref{pvdev-9-th1}}\label{sec3}

\bprop\label{pvdev-9-pro1} For $f\in {\mathcal F}_\alpha(\lambda)$ and
$\lambda\in\mathbb{D}$ we have
\be\label{pvdev-9-eq8}
\left|\frac{f''(z)}{f'(z)}+\alpha f'(z) -c(z,\lambda)\right|\leq r(z,\lambda),
\quad z\in\mathbb{D},
\ee
where
\begin{eqnarray*}
c(z,\lambda)& = & \frac{2[\lambda(1-|z|^2)+\overline{z}(|z|^2-|\lambda|^2)]}
{(1-|z|^2)(1+|z|^2-2{\rm Re\,}({\lambda} z))}, ~~ \mbox{ and }\\
r(z,\lambda) & = & \frac{2(1-|\lambda|^2)|z|}
{(1-|z|^2)(1+|z|^2-2{\rm Re\,}({\lambda} z))}.
\end{eqnarray*}
For each $z\in\mathbb{D}\setminus\{0\}$, equality holds if and only
if $f=H_{e^{i\theta},\lambda}$ for some $\theta\in\mathbb{R}$.
\eprop \bpf
Set $\gamma=0$, $\beta=0$ in  \cite[Proposition 4.1]{samy-vasudev09d}. Then
\cite[Proposition 4.1]{samy-vasudev09d} takes the following form
\begin{eqnarray}\label{pvdev-9-eq8a}
& & \left|1+\frac{zf''(z)}{f'(z)}+\alpha zf'(z) -\frac{(1+\lambda z)(1-\overline{\lambda}\overline{z})
+|z|^2(\overline{z}-\lambda)(\overline{\lambda}+z)}{(1-|z|^2)(1+|z|^2-2{\rm Re\,}({\lambda} z))}\right|\\
\hspace*{1cm}& \leq & \frac{2(1-|\lambda|^2)|z|^2}
{(1-|z|^2)(1+|z|^2-2{\rm Re\,}({\lambda} z))}\nonumber ,
\quad z\in\mathbb{D}.
\end{eqnarray}
A simplification of (\ref{pvdev-9-eq8a}) gives (\ref{pvdev-9-eq8}).
\epf

The choice $\lambda=0$ gives the following result which may need a
special mention.

\bcor\label{pvdev-9-cor01}
For $f\in{\mathcal F}_\alpha(0)$ we have
\be\label{pvdev-9-eq22}
\left|\frac{f''(z)}{f'(z)}+\alpha f'(z)-\frac{2|z|^2\overline{z}}{1-|z|^4}\right|\leq
\frac{2|z|}{1-|z|^4}, \quad z\in\mathbb{D}.
\ee
For each
$z\in\mathbb{D}\setminus\{0\}$, equality holds if and only if
$f=H_{e^{i\theta},0}$ for some $\theta\in\mathbb{R}$.
\ecor

If $f\in{\mathcal F}_\alpha(0)$, then by (\ref{pvdev-9-eq22}) we obtain
$$(1-|z|^2)\left|\frac{f''(z)}{f'(z)}+\alpha f'(z)\right| \leq 2|z|.
$$

\bcor\label{pvdev-9-cor1}
Let $\gamma \colon z(t)$, $0\leq t\leq 1$ be
a $C^1$-curve in $\mathbb{D}$ with $z(0)=0$ and $z(1)=z_0$. Then
we have
$$V (z_0,\lambda)\subset\left\{w\in\mathbb{C}\colon\left |w-
C(\lambda, \gamma) \right |\leq R(\lambda, \gamma)\right\},
$$
where
$$C(\lambda, \gamma)=\int_0^1 c(z(t),\lambda)z'(t)\,dt ~\mbox{ and }~ R(\lambda,
\gamma)=\int_0^1 r(z(t),\lambda)|z'(t)|\,dt.
$$
\ecor
\bpf Since the proof of Corollary \ref{pvdev-9-cor1} follows from \cite{samy-vasudev-vuorinen1}, we omit the
details.
\epf


For the proof of our next result, we need the following lemma.

\blem\label{pvdev-9-lem01}\cite{samy-vasudev-vuorinen1}
For $\theta\in\mathbb{R}$ and $\lambda\in\mathbb{D}$ the function
$$G(z)=\int_0^z \frac{ e^{i\theta}\zeta }
{\{1+(\overline{\lambda}e^{i\theta}-\lambda)\zeta-e^{i\theta}{\zeta}^2\}^2}\,
d\zeta, \quad z\in\mathbb{D},
$$
has a double zero at the origin and no zeros elsewhere in
$\mathbb{D}$. Furthermore there exists a starlike univalent
function $G_0$ in $\mathbb{D}$ such that
$G=e^{i\theta}G^2_0$ and $G_0(0)= G'_0(0)-1=0$.
\elem


\bprop
Let $z_0\in\mathbb{D}\setminus \{0\}$. Then for
$\theta\in(-\pi,\pi]$ we have
$\log H'_{e^{i\theta},\lambda}(z_0)+\alpha H_{e^{i\theta},\lambda}(z_0)
\in\partial V (z_0,\lambda)$.
Furthermore, if  $\log f'(z_0)+\alpha f(z_0)= \log H'_{e^{i\theta},\lambda}(z_0)
+\alpha H_{e^{i\theta},\lambda}(z_0)$ for some $f\in {\mathcal F}_\alpha(\lambda)$
and $\theta\in(-\pi,\pi]$, then $f= H_{e^{i\theta},\lambda}$.
\label{pvdev-9-pro2}
\eprop\bpf
For a proof we refer to \cite[Proposition 4.11]{samy-vasudev09d} with
\vspace{6pt}

\hfill$\displaystyle P(z)=1+\frac{zf''(z)}{f'(z)}+\alpha zf'(z).
$\hfill
\epf

\noindent{\textbf{Proof of Theorem \ref{pvdev-9-th1}}}
Although the proof of Theorem \ref{pvdev-9-th1} is similar to that of the main theorem in
\cite[Theorem 5.1]{samy-vasudev09d},
for the shake of completeness we include the proof. We need to prove that the closed curve
\be\label{pvdev-9-curve}
(-\pi,\pi]\ni \theta \mapsto \log H'_{e^{i\theta},\lambda}(z_0)+
\alpha H_{e^{i\theta},\lambda}(z_0)
\ee
is simple. Suppose that
$$\log H'_{e^{i\theta_1},\lambda}(z_0)+\alpha H_{e^{i\theta_1},\lambda}(z_0)
=\log H'_{e^{i\theta_2},\lambda}(z_0)+\alpha H_{e^{i\theta_2},\lambda}(z_0)
$$
for some $\theta_1,\theta_2\in(-\pi,\pi]$ with $\theta_1\neq\theta_2$.
Then, from Proposition \ref{pvdev-9-pro2}, we have
$$H_{e^{i\theta_1},\lambda}= H_{e^{i\theta_2},\lambda}.
$$
From (\ref{pvdev-9-eq5})  we have
$$
\tau\left(\frac{\omega_{H_{e^{i\theta},\lambda}}}{z},\lambda\right)
=\frac{(1-\overline{\lambda})e^{i\theta} z+\lambda-\overline{\lambda}}
{1-{\lambda}^2-(\lambda-\overline{\lambda})e^{i\theta} z}, ~\tau(z,\lambda) =  \frac{z-\overline{\lambda}}{1-\lambda z}.
$$
That is
$$
\frac{(1-\overline{\lambda})e^{i\theta_1} z+\lambda-\overline{\lambda}}
{1-{\lambda}^2-(\lambda-\overline{\lambda})e^{i\theta_1} z}
=\frac{(1-\overline{\lambda})e^{i\theta_2} z+\lambda-\overline{\lambda}}
{1-{\lambda}^2-(\lambda-\overline{\lambda})e^{i\theta_2} z}
$$
and a simple computation yields
$$
e^{i\theta_1}z=e^{i\theta_2} z
$$
which is a contradiction for the choice of $\theta_1$ and $\theta_2$.
Thus, the curve must be simple. Since $V (z_0,\lambda)$ is a compact
convex subset of $\mathbb{C}$ and has nonempty interior, the boundary
$\partial V (z_0,\lambda)$ is a simple closed curve. From Proposition
\ref{pvdev-9-pro2}, the curve $\partial V (z_0,\lambda)$ contains the curve (\ref{pvdev-9-curve}).
Recall the fact that a simple closed curve cannot contain any simple closed curve
other than itself. Thus, the curve $\partial V (z_0,\lambda)$ is given by (\ref{pvdev-9-curve}).
$\hfill\Box$\\

\section{Geometric view of Theorem \ref{pvdev-9-th1}}\label{sec5}

  Using Mathematica 7,  we describe the boundary of the set
$V(z_0, \lambda)$.  Here we give the Mathematica program which is
used to plot the boundary of the set  $V(z_0, \lambda)$.
We refer to \cite{Ruskeepaa} for the basic concepts on Mathematica programming. The short
notations in this program are of the form: ``z0 for $z_0$'', ``lam for $\lambda$''.

{\tt
\begin{verbatim}

Remove["Global`*"];

z0 = Random[] Exp[I Random[Real, {-Pi, Pi}]]
lam = Random[] Exp[I Random[Real, {-Pi, Pi}]]

Q[lam_, the_] := (2 (Exp[I*the]*z +lam))/
(1 + (Conjugate[lam]*Exp[I*the] - lam)*z - Exp[I*the]*z^2)


myf[lam_, the_, z0_] :=NIntegrate[Q[lam, the], {z, 0, z0},
                       PrecisionGoal -> 2]

image = ParametricPlot[With[{q = myf[lam, the, z0]},
        {Re[q], Im[q]}], {the, -Pi, Pi}]

\end{verbatim}
}

\begin{center}
Table 1
\end{center}

\begin{center}
\begin{tabular}{|c|c|c|}
\hline
Figure & $z_0$ & $\lambda$   \\ \hline
1 & 0.0230875+0.00517512i & 0.175557-0.225417i   \\ \hline
2 & 0.147076+0.0913164i & 0.0748874+0.0476965i    \\ \hline
3 & -0.819143-0.551002i    & 0.722765+0.433556i   \\ \hline
4 & 0.757794-0.598957i   & -0.308071-0.32103i     \\ \hline
5 & -0.414782-0.377338i   & 0.196381-0.500501i     \\ \hline
6 & 0.386456-0.316514i    & -0.236285+0.235873i    \\ \hline
7 & 0.419565+0.478471i    & 0.242605+0.097106i    \\ \hline
8 & 0.754872+0.0830025i   & 0.130907+0.931628i \\ \hline
\end{tabular}
\end{center}

\begin{figure}
\begin{minipage}{0.45\linewidth}
\centering
\includegraphics[width=5.5cm]{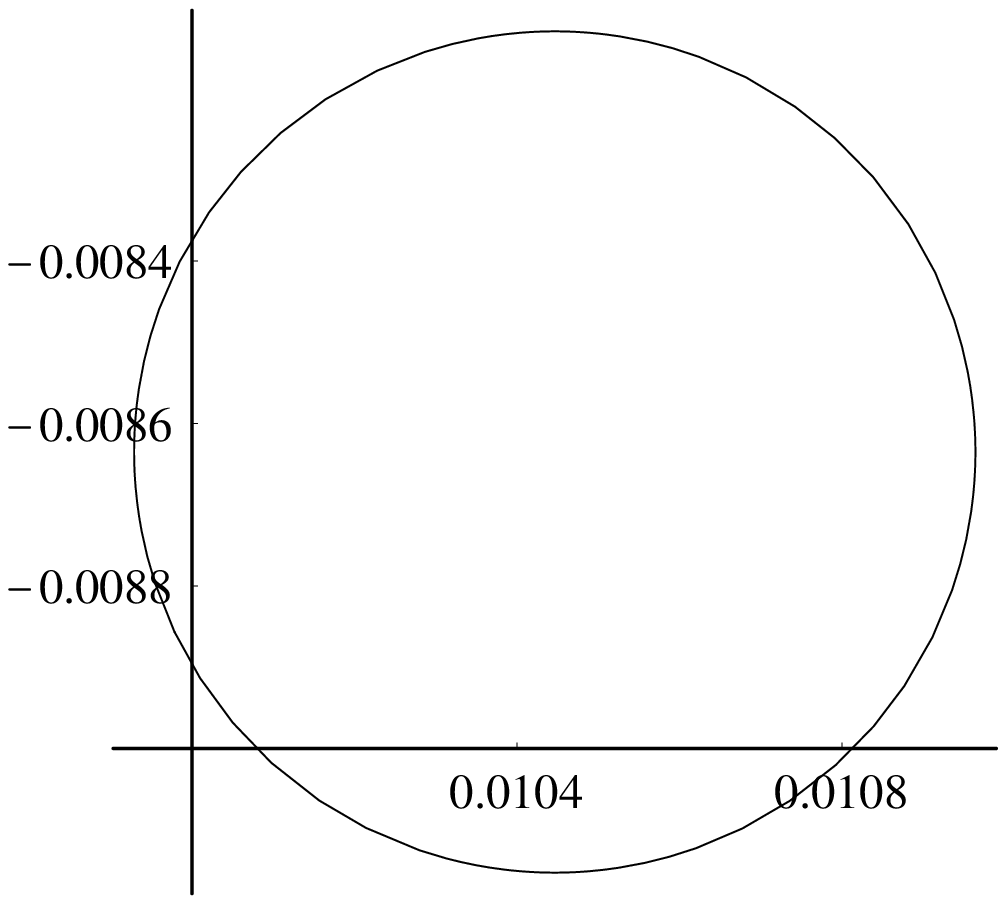}
\caption{}
\end{minipage}
\begin{minipage}{0.5\linewidth}
\centering
\includegraphics[width=5cm]{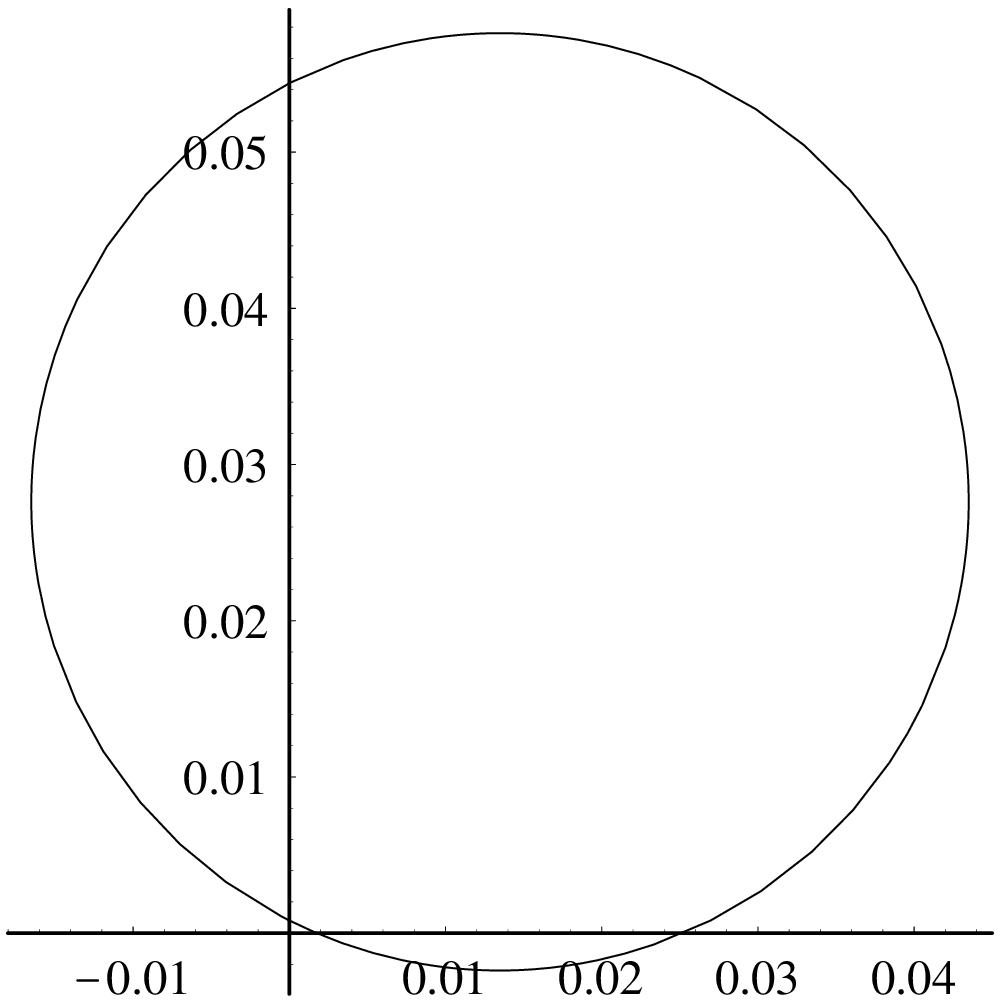}
\caption{}
\end{minipage}
\end{figure}

\begin{figure}
\begin{minipage}{0.45\linewidth}
\centering
\includegraphics[width=6cm]{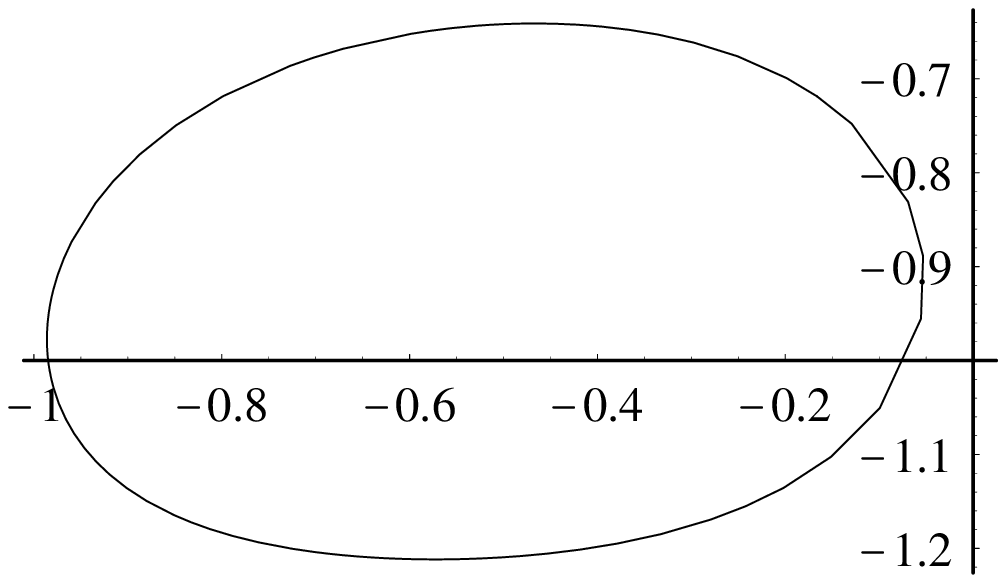}
\caption{}
\end{minipage}
\begin{minipage}{0.5\linewidth}
\centering
\includegraphics[width=5.5cm]{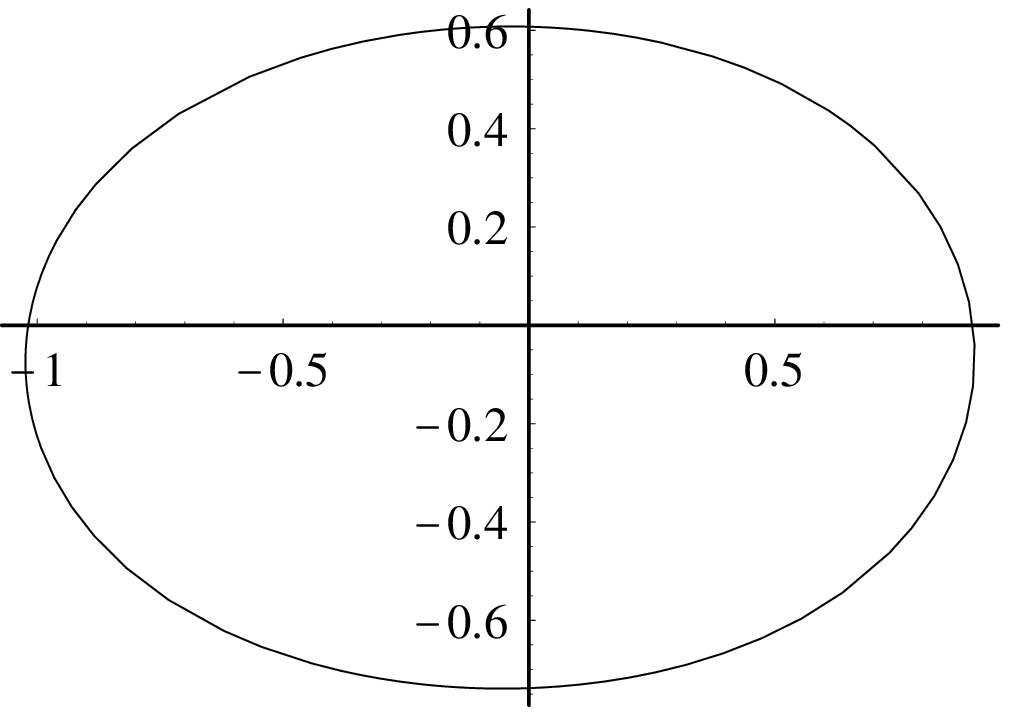}
\caption{}
\end{minipage}
\end{figure}

\begin{figure}
\begin{minipage}{0.45\linewidth}
\centering
\includegraphics[width=5.5cm]{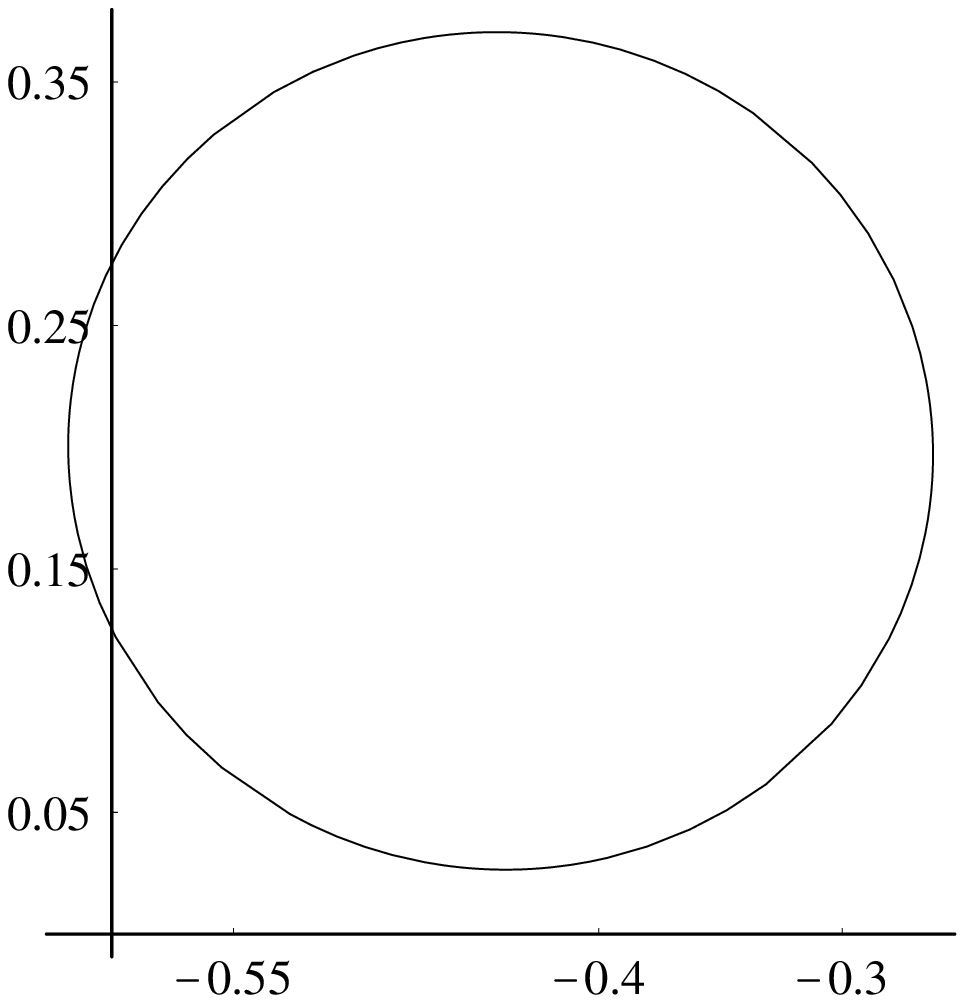}
\caption{}
\end{minipage}
\begin{minipage}{0.5\linewidth}
\centering
\includegraphics[width=5cm]{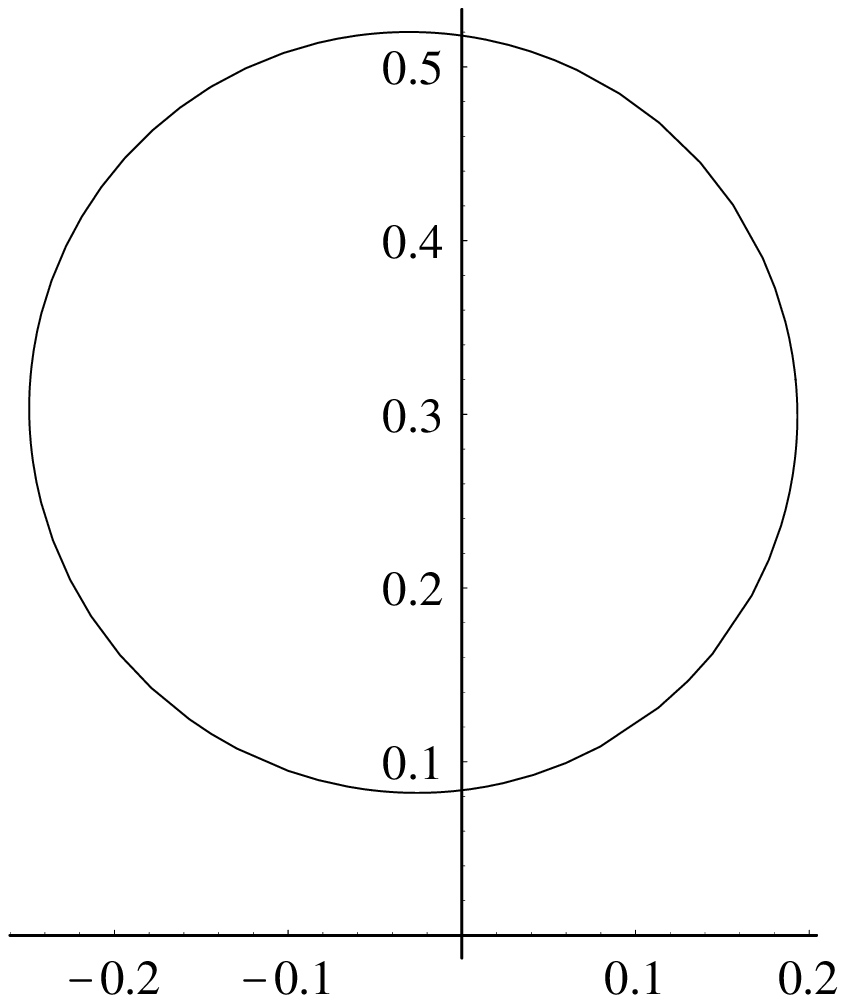}
\caption{}
\end{minipage}
\end{figure}

\begin{figure}
\begin{minipage}{0.45\linewidth}
\centering
\includegraphics[width=5.5cm]{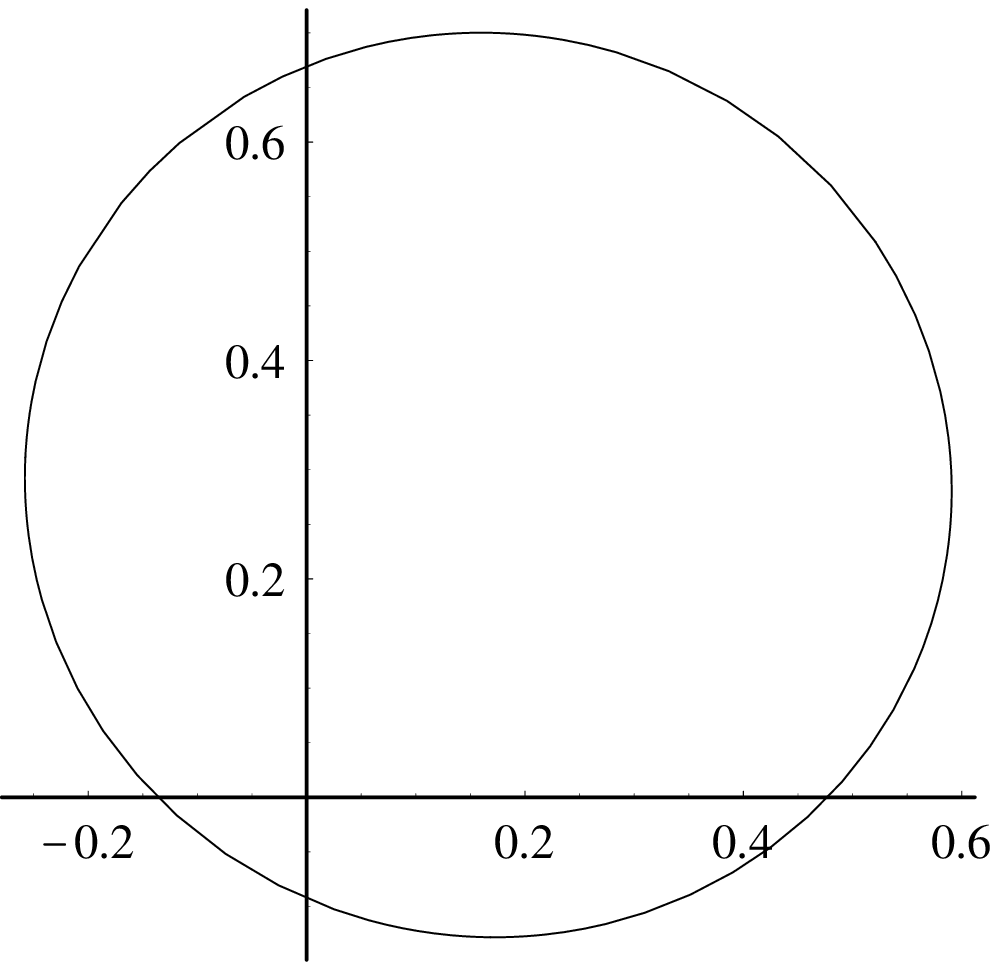}
\caption{}
\end{minipage}
\begin{minipage}{0.5\linewidth}
\centering
\includegraphics[width=5cm]{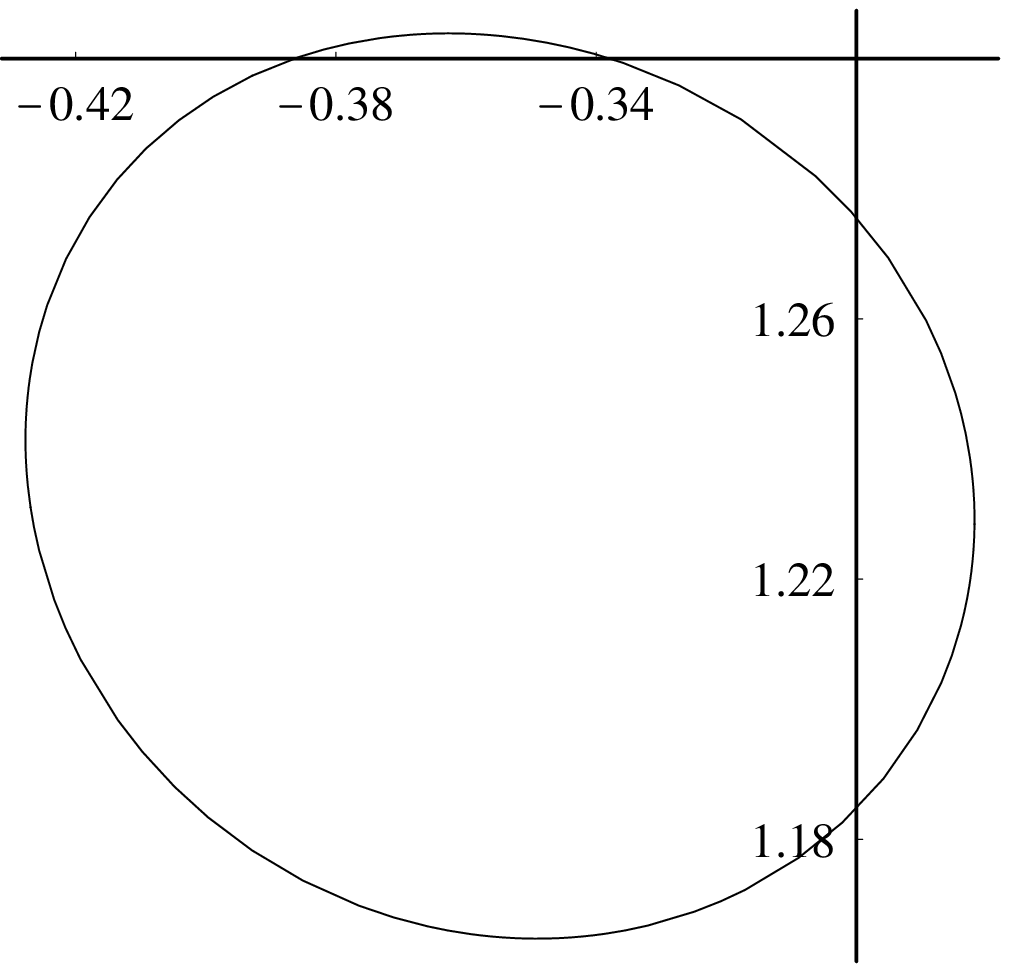}
\caption{}
\end{minipage}
\end{figure}

The Figures from $1$ to $8$ give the geometric view of the region of variability
$V(z_0,\lambda)$ for some sets of parameters  $z_0\in\mathbb{D}$ and
$\lambda\in\mathbb{D}$. We observe from Figures $1$ and $2$ that the regions of variability
$V(z_0,\lambda)$ are very small for some sets of parameters $z_0$ and $\lambda$ whereas
Figures $7$ and $8$ show that the regions of variability $V(z_0,\lambda)$ are
relatively large for some particular sets of parameters.
With the help of Mathematica 7,  we have drawn the curves
$\partial V(z_0, \lambda)$ for various values of $z_0$ and $\lambda$ and observed
that the regions of variability for exponentially convex functions are of the
small size. The above pictures are evident to Proposition {\ref{pvdev-9-pro3}} that
the regions bounded by the curves $\partial V(z_0, \lambda)$ are compact and convex subsets of $\mathbb{C}$.

\section{Open problems}\label{sec7}
\bee
\item  For $f\in\mathcal{E(\alpha)}$, what are the sharp lower and upper bounds of $|f(z)|$ and
$|f'(z)|$ for $z\in\mathbb{D}$?
\item  Let $f\in\mathcal{E(\alpha)}$ and be given by $f(z)=z+\sum_{n=0}^{\infty} a_n z^n$.
Then what are the sharp coefficient bounds for $|a_n|$ for $n\geq 2$?
\eee


\begin{thebibliography}{99}

\bibitem{Arango-Mejia-Rusch-97} {\sc J.H. Arango, D. Mejia, and St. Ruscheweyh,}
Exponentially convex univalent functions,
\textit{Complex Var. Elliptic Equ. } {\bf 33}(1)(1997), 33--50.


\bibitem {Du} {\sc P. L. Duren,} Univalent Functions
(Grundlehren der mathematischen Wissenschaften 259, New York,
Berlin, Heidelberg, Tokyo), Springer-Verlag, 1983.


\bibitem{Goodman's_book} {\sc  A.W. Goodman,}
Univalent Functions, Vols. I and II. Mariner Publishing Co. Tampa,
Florida, 1983.


\bibitem {Pommerenke-book2} {\sc Ch. Pommerenke,} Univalent Functions,
Vandenhoeck and Ruprecht, G\"ottingen, 1975.


\bibitem{samy-book3}{\sc S. Ponnusamy,} Foundations of Complex Analysis,
Alpha Science International Publishers, UK, 2005.

\bibitem{samy-herb} {\sc S. Ponnusamy and H. Silverman,}
Complex Variables with Applications, Birkh{\"{a}}user, Boston, 2006.

\bibitem{samy-vasudev1} {\sc S. Ponnusamy and  A. Vasudevarao,}
Region of variability of two subclasses of univalent functions,
\textit{J. Math. Anal. Appl.} {\bf 332}(2)(2007), 1323--1334.

\bibitem{samy-vasudev09d} {\sc S. Ponnusamy and  A. Vasudevarao,}
Region of variability for functions with positive real part,
\textit{Ann. Polon. Math.} (To appear) 21 pp.

\bibitem{samy-vasudev-vuorinen1}{\sc S. Ponnusamy and A. Vasudevarao, and M. Vuorinen,}
Region of variability for spirallike functions with respect to a boundary point,
\textit{Colloq. Math. }{\bf 116}(1)(2009), 31--46.

\bibitem{samy-vasudev-vuorinen2}{\sc S. Ponnusamy and A. Vasudevarao, and M. Vuorinen},
Region of variability for certain classes of univalent functions
satisfying differential inequalities,
\textit{Complex Var. Elliptic Equ. }{\bf 54}(10)(2009), 899--922.



\bibitem{samy-vasudev-yan4}{\sc S. Ponnusamy, A. Vasudevarao, and H. Yanagihara,}
Region of variability for close-to-convex functions-II,
\textit{Appl. Math. Comput. } {\bf 215}(3)(2009), 901--915.

\bibitem{Ruskeepaa} {\sc H.~Ruskeep\"{a}\"{a},}
Mathematica Navigator: Mathematics, Statistics, and Graphics, $2$nd Ed.,
Elsevier Academic Press, Burlington, MA, 2004.



\bibitem{Yanagihara2} {\sc H. Yanagihara,}
Regions of variability for convex functions,
\textit{Math. Nachr. }{\bf 279}(2006), 1723--1730.

\end{thebibliography}
\end{document}